# AN OVERVIEW OF NUMERICAL AND ANALYTICAL METHODS FOR SOLVING ORDINARY DIFFERENTIAL EQUATIONS

## BYAKATONDA DENIS

## UNIVERSITY OF KISUBI

# Table of Contents

# ABSTRACT

Differential Equations are among the most important Mathematical tools used in creating models in the science, engineering, economics, mathematics, physics, aeronautics, astronomy, dynamics, biology, chemistry, medicine, environmental sciences, social sciences, banking and many other areas **[7].** A differential equation that has only one independent variable is called an Ordinary Differential Equation (ODE), and all derivatives in it are taken with respect to that variable. Most often, the variable is time, t; although, I will use x in this paper as the independent variable. The differential equation where the unknown function depends on two or more variables is referred to as Partial Differential Equations (PDE). Ordinary differential equations can be solved by a variety of methods, analytical and numerical. Although there are many analytic methods for finding the solution of differential equations, there exist quite a number of differential equations that cannot be solved analytically **[8].** This means that the solution cannot be expressed as the sum of a finite number of elementary functions (polynomials, exponentials, trigonometric, and hyperbolic functions). For simple differential equations, it is possible to find closed form solutions **[9].** But many differential equations arising in applications are so complicated that it is sometimes impractical to have solution formulas; or at least if a solution formula is available, it may involve integrals that can be calculated only by using a numerical quadrature formula. In either case, numerical methods provide a powerful alternative tool for solving the differential equations under the prescribed initial condition or conditions **[9].** In this paper, I present the basic and commonly used numerical and analytical methods of solving ordinary differential equations.

# ORGANIZATION OF THE STUDY

Chapter 1 is made up of introduction, which comprises the background of the study, statement of the problem, objectives of the study, justification and the content scope. Chapter 2 highlights on review of literature of ideas of different authors whose findings have been defined in relation to the topic under study. Chapter 3 focuses on methodological review in the light of numerical and analytical methods that are relevant to solving ordinary differential equations, a real world example is also taken as a case study to demonstrate the use of the discussed methods and examine their efficiency. Chapter 4 deals with data analysis and explanation of results obtained from the case studies. In the same way, chapter five consists of a summary, conclusion and recommendations.

The project report however ends with references and appendices in support to the researcher's investigation.

# LIST OF ABBREVIATIONS

ODE        Ordinary Differential Equation

DE        Differential Equation

PDE        Partial Differential Equation

IF        Integrating Factor method

DDEs        Delay Differential Equations

SDEs        Stochastic Differential Equations

LPDEs        Linear elliptic Partial Differential Equations

IVP        Initial Value Problem

BVP        Boundary Value Problem

LHS        Left Hand Side (of an expression)

RHS        Right Hand Side (of an expression)

# LIST OF SYMBOLS

$\mathbb{C}$     The set of complex numbers

$\mathfrak{R}$     Real part of a complex number

$\mathbb{R}$     The set of real numbers

$\mathbb{R}^d$     The d-dimensional vector-space of the real d-tuple

$\mathcal{U}$     The exact solution of an ODE or IVP

$\mathcal{U}_k$     The exact solution at time step $t_k$

$y_k$     The discrete solution at time step $t_k$

$|x|$     The absolute value of a real number, the modulus of a complex number, or the Euclidean norm in $\mathbb{R}^d$ or $\mathbb{C}^d$, depending on its argument

# LIST OF FIGURES

# LIST OF TABLES

# CHAPTER ONE

# INTRODUCTION

## 1.0 Introduction

This chapter focuses on the background of the study, statement of the problem, objectives and purpose of the study, significance of the study and the scope of study. It also highlights the key preliminaries in the study and outlines how the study will be carried out.

## 1.1 Background of the study

A differential equation is an equation that relates a function to its derivative(s). The unknown is the function. A differential equation is said to be ordinary if the function is uni-variate and more precisely if its domain is a connected subset of $\mathbb{R}$. **[27]**

Ordinary differential equations arise in many different contexts. These different contexts include fundamental laws of physics, mechanics, electricity, thermodynamics, and also population and growth modeling **[28]**.

Ordinary differential equations (ODEs) arise in many contexts of mathematics and social and natural sciences, for example, in physics, the Legendre DE **[25]**, which is a self-adjoint ODE, arises in the solutions of Hydrogen atom wave functions and angular momentum in single-particle quantum mechanics. Their solutions form the polar angle part of the spherical harmonics basis for the multi pole expansion, which is used in both electromagnetic and gravitational statics. In engineering, for example, many difficult problems in the field of static and dynamic mechanics can be solved by computing the solutions self-adjoint Bessel equations **[26].**

Mathematical descriptions of change use differentials and derivatives. Various differentials, derivatives, and functions become related via equations, such that a differential equation is a result that describes dynamically changing phenomena, evolution, and variation. Often, quantities are defined as the rate of change of other quantities (for example, derivatives of displacement with respect to time), or gradients of quantities, which is how they enter differential equations.

It is a common truth that Differential Equations are among the most important Mathematical tools used in producing models in the engineering, mathematics, physics, aeronautics, elasticity, astronomy, dynamics, biology, chemistry, medicine, environmental sciences, social sciences, banking and many other areas **[29]**.

During the past few decades, the development of nonlinear analysis, dynamical systems and their applications to science and engineering has stimulated renewed enthusiasm for the theory of Ordinary Differential Equations (ODE).

Differential equations and mathematical modeling can be used to study a wide range of social issues. Among the topics that have a natural fit with the mathematics in a course on ordinary differential equations are all aspects of population problems: growth of population, over-population, carrying capacity of an ecosystem, the effect of harvesting, such as hunting or fishing, on a population and how over-harvesting can lead to species extinction, interactions between multiple species populations, such as predator-prey, cooperative and competitive species **[7]**.

Ordinary differential equations are ubiquitous in science and engineering: in geometry and mechanics, in chemical re-action kinetics, molecular dynamics, electronic circuits, population dynamics, and many more application areas. They also arise, after semi-discretization in space, in the numerical treatment of time-dependent partial differential equations, which are even more impressively omnipresent in our technologically developed and financially controlled world.

The most common specific fields that require modeling in terms of differential equations include geometry and analytical mechanics. Scientific fields include much of physics and astronomy (celestial mechanics), meteorology (weather modeling), chemistry (reaction rates), biology (infectious diseases, genetic variation), ecology and population modeling (population competition), economics (stock trends, interest rates and the market equilibrium price changes).

Finding solutions to differential equations allows us to make reasonable prediction about most natural phenomena, the problem here is what method(s) do we have to use to solve an ODE?, which is the simplest method? And how accurate is the method? All these questions have to be taken into account before solving an ODE and this research will help to answer these questions by classifying the different available methods according to different case usages and comparing their relative efficiency so that a novice mathematician or other researchers find it simple to choose a method to use when they encounter a Differential Equation.

## 1.2 Preliminaries

A function F of $x$ and $y$ and derivatives of $y$ such that

$$F(x, y, y', y'', \ldots\ldots y^{(n-1)}, y^{(n)}) = 0,$$

is called an ordinary differential equation of *order n* for the unknown function *y*.

The ordinary differential equation is said to be *linear* if it can be written as a linear combination of the derivatives of *y*, i.e,

$$y^{(n)} = \sum_{i=0}^{n-1} a_i(x) y^{(i)} + r(x),$$

where $a_i(x)$ and $r(x)$ are continuous functions of *x*. The above equation is *homogeneous* if $r(x) = 0$, *otherwise it is non-homogeneous.*

Among ordinary differential equations, linear differential equations play a prominent role for several reasons. Most elementary and special functions that are encountered in physics and applied mathematics are solutions of linear differential equations. When physical phenomena are modeled with non-linear equations, they are generally approximated by linear differential equations for an easier solution. The few non-linear ODEs that can be solved explicitly are generally solved by transforming the equation into an equivalent linear ODE **[2]**

Linear differential equations, which have solutions that can be added and multiplied by coefficients, are well-defined and understood, and exact closed-form solutions are obtained. By contrast, ODEs that lack additive solutions are nonlinear, and solving them is far more intricate, as one can rarely represent them by elementary functions in closed form: Instead, exact and analytic solutions of ODEs are in series or integral form. Numerical methods applied by hand or by computer, may approximate solutions of ODEs and perhaps yield useful information, often sufficing in the absence of exact, analytic solutions.

## 1.4 Statement of the problem

Many studies have been devoted in order to find the solution of ordinary differential equations. In the case where the equation is linear, it is not a major problem since it can be solved by analytical methods. Unfortunately, most of the interesting differential equations that are got after modeling real world problems are non-linear and it causes a major problem to solve that equation by analytical methods. Thus, numerical methods have been developed and have proved really helpful to solve those ordinary differential equations. Furthermore, there is much computer software that has been developed to help the user to solve those equations.

Inasmuch as Ordinary differential equations frequently occur as mathematical models in many branches of science, engineering and economy. It is Unfortunate that these equations seldom have solutions that can be expressed in closed form, so it is common to seek approximate solutions by means of numerical methods; This research will investigate and compare the various analytical and numerical methods of solving ordinary differential equations within limits of low error bounds including computer-enhanced techniques to solve problems that would otherwise be impossible to solve or take a lot of time.

## 1.5 Objectives of the study

### 1.5.1 General objective

- To solve differential equations using analytic and numerical methods and plot the graphs of the solutions to obtain qualitative information about the problem.

### 1.5.2 Specific objectives

- To perform a comparative study on analytical and numerical solutions of Initial Value Problems (IVP) for Ordinary Differential Equations (ODE).
- To compute the error between the numerical solutions obtained by the Euler method, Runge-kutta methods and the analytical solutions.

## 1.6 Significance of the study

Differential equations are commonly used for mathematical modeling in science and engineering. Many problems of mathematical physics can be stated in the form of differential equations. These equations also occur as reformulations of other mathematical problems such as ordinary differential equations and partial differential equations. In most real life situations, the differential equation that models the problem is too complicated to solve exactly, and one of two approaches is taken to approximate the solution. The first approach is to simplify the differential equation to one that can be solved exactly and then use the solution of the simplified equation to approximate the solution to the original equation. The other approach uses numerical methods for approximating the solution of original problem. This is the approach that is most commonly taken since the approximation methods give more accurate results and realistic error information. Numerical methods are generally used for solving mathematical problems that are formulated in science and engineering where it is difficult or even impossible to obtain exact solutions. Only a limited number of differential equations can be solved analytically. **[29]**

There are many analytical methods for finding the solution of ordinary differential equations but even then there exist a large number of ordinary differential equations whose solutions cannot be obtained in closed form by using the well-known analytical methods, where we have to use the numerical methods to get the approximate solution of a differential equation under the prescribed initial condition or conditions. There are many types of practical numerical methods for solving initial value problems for ordinary differential equations. In this research I present, analyze and compare two standard numerical methods Euler and Runge-Kutta for solving initial value problems of ordinary differential equations.

## 1.7 Justification

The study will afford mathematicians the opportunity to be aware of some key numerical and analytical methods for solving ordinary differential equations. Information gathered from the results would educate mathematicians and other researchers about the stability, efficiency and case usage of these methods and which one of these methods is much more preferred than the others.

This will help fill the gap in the research carried out in other universities in Uganda and all around the world. In addition, it could pave the way for more comprehensive research on the comparison of these methods in relation to some specified complex functions which are very significant in drawing conclusions to research works.

Researchers will be alerted of these important methods which play a very important role in obtaining accurate results when used to compute stiff differential equations.

The study will equally be helpful to natural scientists (chemists, physicists and biologists) and engineers to understand which of these numerical and analytical schemes would be suitable to solve a modeled problem if it is an Ordinary Differential Equation (ODE).

## 1.8 Scope of the study

Differential equations are of many categories and each category requires its own techniques and methods of solving it. The two broad categories of differential equations are the ordinary differential equations and the partial differential equations. Other subcategories include Delay Differential Equations (DDEs), Stochastic Differential Equations (SDEs), linear elliptic Partial Differential Equations (LPDEs) and others.

However, this research will only investigate the analytical and numerical methods of solving ordinary differential equations, linear and non-linear, homogeneous and non-homogeneous. More concentration will be put in comparing the solutions obtained from Euler method and the Runge-Kutta methods with the analytical solutions, the superiority of numerical methods will also be compared in terms of rate of convergence for a given step length.

# CHAPTER 2

# LITERATURE REVIEW

## 2.0 Introduction

This chapter begins with the history of the evolution of ordinary differential equations and the search for better methods of finding their solutions, it summarizes the different numerical and analytical approaches used in solving ordinary differential equations as put forward by several writers and researchers who have investigated this study before or have conducted similar research.

## 2.1 The necessity of a solution to a differential equation

Differential equations can describe nearly all systems undergoing change. They are ubiquitous is science and engineering as well as economics, social science, biology, business, health care, etc. Many researchers and mathematicians have studied the nature of Differential Equations and many complicated systems that can be described quite precisely with mathematical expressions.

Whenever we meet or form a differential equation, a number of questions arise right away:

1. Does the equation have a solution?
2. If it does, is the solution unique?
3. What is the solution?
4. Is there a systematic way to solve such an equation?

Solutions of ordinary differential equations (ODEs) are in general possible by different methods **[8]**. The main methods of solving ordinary differential equations are analytical and numerical; all other approaches are subsets of these.

## 2.2 History of Ordinary Differential Equations

### 2.2.1 Overview

The attempt to solve physical problems led gradually to mathematical models involving an equation in which a function and its derivatives play important roles. However, the theoretical development of this new branch of mathematics - Ordinary Differential Equations - has its origins rooted in a small number of mathematical problems. These problems and their solutions led to an independent discipline with the solution of such equations an end in itself.

In *circa* 1671, English physicist Isaac Newton wrote his then-unpublished *The Method of Fluxions and Infinite Series* (published in 1736), in which he classified first order differential equations, known to him as *fluxional equations*, into three classes **[1],** as follows (using modern notation):

| Ordinary differential equations | | Partial differential equations |
|---|---|---|
| Class 1 | Class 2 | Class 3 |
| $\dfrac{dy}{dx} = f(x)$ | $\dfrac{dy}{dx} = f(x,y)$ | $x\dfrac{\partial u}{\partial x} + y\dfrac{\partial u}{\partial y} = u$ |

The first two classes contain only ordinary derivatives of one or more dependent variables, with respect to a single independent variable, and are known today as "*ordinary differential equations*"; the third classes involves partial derivatives of one dependent variable and today are called "*partial differential equations*".

The study of "differential equations", according to British mathematician Edward Ince, and other mathematics historians is said to have begun in 1675, when German mathematician Gottfried Leibniz wrote the following equation **[2]**

$$\int x\,dx = \frac{1}{2}x^2.$$

In 1676, Newton solved his first differential equation. That same year, Leibniz introduced the term "differential equations" (*aequatio differentialis,* Latin) or to denote a relationship between the differentials *dx* and *dy* of two variables *x* and *y* **[6]**

In 1693, Leibniz solved his first differential equation and that same year Newton published the results of previous differential equation solution methods—a year that is said to mark the inception for the differential equations as a distinct field in mathematics.

Swiss mathematicians, brothers Jacob Bernoulli (1654-1705) and Johann Bernoulli (1667-1748), in Basel, Switzerland, were among the first interpreters of Leibniz' version of differential calculus. They were both critical of Newton's theories and maintained that Newton's theory of fluxions was plagiarized from Leibniz' original theories, and went to great lengths, using differential calculus, to disprove Newton's *Principia*, on account that the brothers could not accept the theory, which Newton had proven, that the earth and the planets rotate around the sun in elliptical orbits. **[3]**

The first book on the subject of differential equations, supposedly, was Italian mathematician Gabriele Manfredi's 1707 *On the Construction of First-degree Differential Equations*, written between 1701 and 1704, published in Latin. **[4]** The book was largely based or themed on the views of the Leibniz and the Bernoulli brothers. Most of the publications on differential equations and partial differential equations, in the years to follow, in the 18th century, seemed to expand on the version developed by Leibniz, a methodology, employed by those as Leonhard Euler, Daniel Bernoulli, Joseph Lagrange, and Pierre Laplace.

In 1739, Swiss mathematician Leonhard Euler began using the integrating factor as an aid to derive differential equations that were integrable in finite form **[5]**

The circa 1828 work of English physical mathematician George Green seems to have something to do with defining a test for an "integrable" or conservative field of force (or somehow has connection to thermodynamics via William Thomson); such as in terms of the later 1871 restylized "curl" notation (test of integrability) of James Maxwell (or possibly the earlier work of Peter Tait). **[10]** In circa 1839, Green stated:

*"If all the internal forces exerted be multiplied by the elements of their respective directions, the total sum for any assigned portion of the mass will always be the exact differential of some function."*

The strain-energy potential function of Green is said to of the same theme as Willard Gibbs thermodynamics potentials and Hermann Helmholtz free energy.**[12]** The use of the both terms "exact differential" and "complete differential" were in common use at least as early as 1841 **[11]**.

From 1850 to 1875, German physicist Rudolf Clausius revolutionized physical science (chemistry, physics, and mechanics) when he transformed the failing notion of French chemist Antoine Lavoisier's "caloric particle model of heat"—in which a single differential unit or quantity of heat was considered to be an small fluid-like particle (smaller in size than an atom) that was indestructible and said to be located in the interstices of bodies (in the space between the atoms) in various amounts, dependent upon the volume of the given body (more in the body for large volumes; less for smaller volumes) according to Boerhaave's law—into that of a quantity of heat dQ defined as the product of the absolute temperature T of a body and the "exact differential" quantity entropy $dS$, such that $dQ = TdS$, and the physical-mathematical function $dQ/T$ is an extensive exact differential quantity state function. This is probably the most complicated mathematical formalisms in all of human knowledge.

Clausius began to introduce some of the mathematical background to this notion of the "exact differential model of heat" in his 1858 article "On the Treatment of Differential Equations which are Not Directly Integrable", in which he introduced the now-infamous "condition for an exact differential" to justify his claim that $1/T$ is the integrating factor (T being the integrating denominator) of the inexact differential function dQ, which makes the resulting function dQ/T an exact differential.

### 2.2.2 Evolution of the methods of solving ordinary Differential Equations as seen in literature

The search for general methods of integrating differential equations began when Isaac Newton (1642-1727) classified first order differential equations into three classes **[1]**. Newton would express the right side of the equation in powers of the dependent variables and assumed as a solution an infinite series. The coefficients of the infinite series were then determined **[13]**.

Even though Newton noted that the constant coefficient could be chosen in an arbitrary manner and concluded that the equation possessed an infinite number of particular solutions, it wasn't until the middle of the 18th century that the full significance of this fact, i.e., that the general solution of a first order equation depends upon an arbitrary constant, was realized.

In 1692 James Bernoulli made known the method of integrating the homogeneous differential equation of the first order, and not long afterwards reduced to quadratures the problem of integrating a linear equation of the first order **[14]**.

The original discoveries of practically all known elementary methods of solving differential equations of the first-order took place during the Bernoulli dynasty.

The Bernoullis' were a Swiss family of scholars whose contributions to differential equations spanned the late seventeenth and the eighteenth century. Nikolaus Bernoulli I (1623-1708) was the progenitor of this celebrated family of mathematicians. James I, John I, and Daniel I are the best-known members of the Bernoulli family who made many contributions to this new field of differential equations.

In 1691 the inverse problem of tangents led Leibniz to the implicit discovery of the method of separation of variables **[15]**. However, it was John Bernoulli, in a letter to Leibniz, dated May 9, 1694, that gave us the explicit process and the term, *seperatio indeterminatarum* or separation of variables **[16]**.

But even then, in one, but important, case of:

$$x\, dy - y\, dx = 0.$$

This process failed, because it led to $\frac{1}{y} dy = \frac{1}{x} dx$ and $\frac{1}{x} dx$ had not yet been integrated.

In 1696, John Bernoulli, a student, rival, and equal of his older brother James, gave a main impetus to the study of differential equations through posing his famous brachistochrone problem of finding the equation of the path down which a particle will fall from one point to another in the shortest time [17].

The equation

$$\frac{dy}{dx} + P(x)y = Q(x)y^n,$$

known today as the Bernoulli equation, was proposed for solution by James Bernoulli in December, 1695 [16]. The following year Leibniz solved the equation by making substitutions and simplifying to a linear equation, similar to the method employed today [18].

In the early years of the eighteenth century a number of problems led to differential equations of the second and third order. In 1701 James Bernoulli published the solution to the *isoperimetric problem* - a problem in which it is required to make one integral a maximum or minimum while keeping constant the integral of a second given function - thus resulting in a differential equation of the third order. [19, 18]. In a letter written to Leibniz, May 20, 1716, John Bernoulli discussed the equation:

$$\frac{d^2y}{dx^2} = \frac{2y}{x^2},$$

where the general solution when written in the form

$$y = \frac{x^2}{a} + \frac{b^2}{3x},$$

involves three cases: When *b* approaches zero the curves are parabolas; when *a* approaches infinity, they are hyperbolas; otherwise, they are of the third order [19].

Jacopo Riccati's (1676-1754) discussion of special cases of curves whose radii of curvature were dependent solely upon the corresponding ordinates resulted in his name being associated with the equation

$$y' = P(x) + Q(x)y + R(x)y^2.$$

Riccati's discussion offered no solutions of his own. It was Daniel Bernoulli who successfully treated the equation commonly known as the *Riccati Equation*. **[20].**

In general, this equation cannot be solved by elementary methods. However, if a particular solution $y_1(x)$ is known, then the general solution has the form:

$$y(x) = y_1(x) + z(x),$$

where z(x) is the general solution of the Bernoulli equation:

$$z' - (Q + 2Ry_1)z = Rz^2.$$

By 1724 Daniel Bernoulli had found the necessary and sufficient conditions for integrating in a finite form the equation:

$$y' + ay^2 = bx^m.$$

Leonhard Euler (1707-1783) provided the next significant development when he posed and solved the problem of *reducing a particular class of second order differential equations to that of first order.* His process of finding a second solution from a known solution consists both of reducing a second order equation to a first order equation and of *finding an integrating factor.* Additionally, Euler confirmed that the ratio of two different integrating factors of a first-order differential equation is a solution of the equation **[21].**

Euler began his treatment of the homogeneous linear differential equation with constant coefficients in a letter he wrote to John Bernoulli on September 15, 1739, published in *Miscellanea Berolinensia*, 1743. Within a year Euler had completed this treatment by successfully dealing with repeated quadratic factors and turned his attention to the non-homogeneous linear equation **[21].**

The method of successive reduction of the order of the equation with the aid of integrating factors led first to equations integrable in finite form. Euler first reduced these equations step by step and then integrated. For those equations which were not integrable in a finite form, Euler used the method of integrating by series.

The famous method of Euler was published in his three-volume work *Institutiones Calculi Integralis* in the years 1768 to 1770, republished in his collected works (Euler, 1913).

Joseph Louis Lagrange (1736-1813), while working on the *problem of determining an integrating factor for the general linear equation,* formalized the concept of the *adjoint equation*. Lagrange not only determined an integrating factor for the general linear equation, but furnished proof of the general solution of a homogeneous linear equation of order *n* **[22].** In addition, Lagrange discovered the method of variation of parameters. **[23].**

Building on Lagrange's work, Jean Le Rond d'Alembert (1717-1783), found the *conditions under which the order of a linear differential equation could be lowered* **[10].** By deriving a method of dealing with the exceptional cases, d'Alembert solved the *problem of linear equations with constant coefficients*, and initiated the study of *linear differential systems* **[24].** In a treatise written in 1747, devoted to vibrating strings, d'Alembert was led to partial differential equations where he did his main work in the field.

The period of initial discovery of general methods of integrating ordinary differential equations ended by 1775, a hundred years after Leibniz inaugurated the integral sign. For many problems the formal methods were not sufficient. Solutions with special properties were required, and thus, criteria guaranteeing the existence of such solutions became increasingly important. Boundary value problems led to ordinary differential equations, such as Bessel's equation, that prompted the study of Laguerre, Legendre, and Hermite polynomials. The study of these and other functions that are solutions of equations of hyper geometric type led in turn to modern numerical methods.

Thus, by 1775, as more and more attention was given to analytical methods and problems of existence, the search for general methods of integrating ordinary differential equations ended momentarily.

The idea of generalizing the Euler method, by allowing for a number of evaluations of the derivative to take place in a step, is generally attributed to Runge (1895). Further contributions were made by Heun (1900) and Kutta (1901). The latter completely characterized the set of Runge–Kutta methods of order 4, and proposed the first methods of order 5. Special methods for second order differential equations were proposed by Nyström (1925), who also contributed to the development of methods for first order equations. It was not until the work of Huˇta (1956, 1957) that sixth order methods were introduced.

Islam, Md.A. (2015), published a journal on the accuracy analysis of Numerical solutions of Initial Value Problems (IVP) for Ordinary Differential Equations (ODE) **[30],** Shampine, L.F. and Watts, H.A. (1971) compared the error estimators for Runge-Kutta Methods **[33]**

Separation of variable method of solving partial differential equation is also called Fourier's method [Renze, John and Weisstein, Eric W. **[36]**" Separation of Variables]. This method is effective because of the fact that if the product of functions of independent variables is a constant, each function must separately be a constant. The first person to use the method was L'Hospital in 1750 [Renze, John and Weisstein, Eric W." Separation of Variables].

Pierre-Simon Laplace, a French Mathematician introduced a special type of integral transform in his research later on it was called as Laplace transformation. Oliver Heaviside, a British physicist, developed Laplace transformation systematically. Among the various integral transform, it is used mostly. The easiness in understanding and simple in applying is the inspiration behind this transformation technique. In numerous problems Laplace transformation is applied to derive the general solution.

From this review, the history of the search for better methods of solving ordinary differential equations has come a long way and quite a number of methods have been established, there's however need to study these methods and compare them in terms of efficiency and classify them to an appropriate usage (especially for people who are simply looking for a better method to solve their ODE, and not necessarily studying the subject in depth), something which the previous researchers on this subject did not take consideration of.

Also from the literature review, several works in numerical solutions of initial value problems using Euler method and Runge-Kutta method have been carried out. Many authors have attempted to solve initial value problems (IVP) to obtain high accuracy with speed by using numerous methods, such as Euler method, midpoint method and Runge-Kutta method, and also some other methods. In **[30]** the author discussed accuracy analysis of numerical solutions of initial value problems (IVP) for ordinary differential equations (ODE), and also in **[31]** the author discussed accurate solutions of initial value problems for ordinary differential equations with fourth-order Runge-kutta method. In **[32]** some numerical methods for solving initial value problems in ordinary differential equations are studied. In **[33]** numerical solutions of initial value problems for ordinary differential equations were studied using various numerical methods.

In this research, Euler method and Runge-Kutta method are applied without any discretization, transformation or restrictive assumptions for solving ordinary differential equations in initial value problems.

According to the researcher's knowledge and from the review, a lot of studies have been done about the analytical methods and the numerical methods of solving ODEs but these topics are usually disjoint in literature. The contribution of this study is to put together different analytical methods and numerical methods for solving ODES.

# CHAPTER 3:

# RESEARCH METHODOLOGY

## 3.0 Introduction

This chapter describes the methods and instruments that were used by the researcher to demonstrate the various numerical and analytical methods for solving ordinary differential equations (ODEs) and show their comparative efficiencies in a case use.

Data analysis, simulation and modeling software such as MATLAB, was used. Specifically, the results presented in chapter 4 were generated using MATLAB

All programs for solving Initial value problems with numerical methods in this study were be written in MATLAB since its easier and specific to scientific computing, it also makes graphical analysis of the solutions much faster.

## 3.1 Some important theorems, definitions and lemmas

Here I introduce a few important definitions, Lemmas and theorems that relate to ordinary differential equations. Most of these will be referred to throughout the text.

*Definition 3.1.1* (normal form of a differential equation):

Considering an $n^{th}$ order linear differential equation is an equation of the form

$$a_0(x)y + a_1(x)y^{(')} + a_2(x)y^{('')} + \cdots + a_n(x)y^{(n)} + b(x) = 0$$

where $a_0(x) \ldots a_n(x)$ and $b(x)$ are arbitrary differentiable functions of x.

Dividing through by $a_n(x)$, the general $n^{th}$ order linear equation can be written in normal form as

$$y^{(n)} + \cdots + P_2(x)y^{('')} + P_1(x)y^{(')} + P_0(x)y + Q(x) = 0$$

where $P_n(x) = \frac{a_i(x)}{a_n(x)}$ .

*Definition 3.1.2 (Solution to a differential equation)*: A real valued function $y(x)$ defined on an interval $I \subset R$, is called a (particular) solution to a differential equation if $y(x)$ is differentiable at any $x \in I$, the point $(x, y(x))$ belongs to $D$ for any $x \in I$ and the identity $y'(x) = f(x, y(x))$ holds for all $x \in I$. The family of all particular solutions of a differential equation is called the **general solution**. The graph of a particular solution is called an integral curve of the equation. Obviously, any integral curve is contained in the domain $D$.

Below are some examples of differential equations,

$$\frac{dy}{dx} = 2x + 1.$$

$$\frac{dy}{dx} - \frac{y}{x} = x^2.$$

$$\frac{d^2y}{dx^2} + y = 0.$$

$$\frac{d^2y}{dx^2} + 3\frac{dy}{dx} + 2y = \sin x.$$

The solutions to each of the above equations involve finding the unknown function $y(x)$ such that it satisfies the equation. Solutions that involve arbitrary constants (which depend on the order of the ODE) are called *general solutions* because they represent the whole family of possible solutions. A solution which satisfies the differential equation but contains no arbitrary constants is called a *particular solution*, such a solution can be got when some initial conditions are given, the so-called *initial value problem* (IVP).

*Definition* 3.1.3 (Wronskian): Let $y_1, y_2 \ldots, y_n$ are $n$ real-valued functions on $I$, which are $(n-1)$ times differentiable on $I$ then their Wronskian is defined by

$$W(x) = \det \begin{pmatrix} y_1 & y_2 & \cdots & y_n \\ y_1' & y_2' & \cdots & y_n' \\ \cdots & \cdots & & \cdots \\ y_1^{(n-1)} & y_2^{(n-1)} & \cdots & y_n^{(n-1)} \end{pmatrix}$$

*Theorem 3.1.1* (Abel's formula): For an nth order linear equation, the wronksian W is given by,

$$W = Ae^{-\int P_1(x)dx} \qquad [37]$$

## 3.2 Analytical methods of solving an ordinary differential equation

Analytical methods as the name suggests involve 'analysis' of the differential equation such that it can be made simple separable or an exact differential be formed on one side (for the case of first order ordinary differential equations), For higher order equations, analysis enables us to reduce the order to one that is easier to solve. Analytical methods vary depending on the order and form of differential equation, inasmuch as they are very useful and produce accurate results, they have a few limitations such as the difficulty in analysis and furthermore, not all differential equations can be solved using these methods.

Here I present the different analytical methods and the kind of ordinary differential equations they can solve.

### 3.2.0   First order ODEs

**(a) Separation of variables**

Let's consider an ODE of the form

$$\frac{dy}{dx} = \frac{f(x)}{g(y)}.$$

Such an equation is said to have **separable variables,** because it can be rearranged to get just $y$ on the left side and just $x$ on the right. This process is called 'separating the variables'.

Multiplying by $g(y)$ gives,

$$g(y)\frac{dy}{dx} = f(x).$$

This is the kind of equation we get when an implicit function is differentiated. If we can find functions $G(y)$ and $F(x)$ such that $G'(y) = g(y)$ and $F'(x) = f(x)$, then the equation can be written as

$$G'(y)\frac{dy}{dx} = F'(y).$$

By analysis, the term on the LHS is $\frac{d}{dx}G(x)$. Thus

$$\frac{d}{dx}G(x) = F'(y),$$

integrating both sides, we obtain

$$G(x) = F(y) + C.$$

### (b) Analysis to find an exact differential

If the terms on either the LHS or RHS of the ordinary differential equation are seen to be forming an exact differential, then the ODE can be solved by integrating both sides of the equation.

Using algebra, any first order equation can be written in the form

$F(x,y)dx + G(x,y)dy = 0$, for some functions $F(x,y), G(x,y)$.

*Definition:* An expression of the form $F(x,y)\,dx + G(x,y)\,dy$ is called a (first-order) differential form. A differential form $F(x,y)\,dx + G(x,y)\,dy$ is called exact if there exists a function g(x, y) such that $dg = F\,dx + G\,dy$.

*Theorem:* If F and G are functions that are continuously differentiable throughout a simply connected region, then $F\,dx + G\,dy$ is exact if and only if $\frac{\partial G}{\partial x} = \frac{\partial F}{\partial y}$ [38].

### (c) Integrating factor method

In case the first order linear ODE is neither simple separable nor exact, then, it can be made exact by multiplying through with a term called an *'integrating factor'*

Consider a first order linear non-homogeneous ODE below, and assume it is inexact

$$a_1(x)\frac{dy}{dx} + a_0(x)y = f(x).$$

Since $a_1(x) \neq 0$; we can divide through by it to obtain

$$\frac{dy}{dx} + \frac{a_0(x)}{a_1(x)}y = \frac{f(x)}{a_1(x)},$$

and then let

$$\frac{a_0(x)}{a_1(x)} = P(x) \text{ and } \frac{f(x)}{a_1(x)} = Q(x).$$

Thus

$$\frac{dy}{dx} + P(x)y = Q(x).$$

A first order linear inexact ODE of the special form $\frac{dy}{dx} + P(x)y = Q(x)$ can be made exact by multiplying through by an integrating factor (IF).

Let the integrating factor be $\mu(x)$ and suppose multiplying $\frac{dy}{dx} + P(x)y = Q(x)$ by $\mu(x)$ makes it exact;

$$\mu(x)\frac{dy}{dx} + \mu(x)P(x)y = \mu(x)Q(x).$$

The aim of this method is to choose $\mu(x)$ so that

$$\mu(x)\frac{dy}{dx} + \mu(x)P(x)y = \frac{d}{dx}[\mu(x)y].$$

Since

$$\frac{d}{dx}[\mu(x)y] = \mu(x)\frac{dy}{dx} + y\frac{d}{dx}\mu(x).$$

Thus $\mu(x)$ must be chosen so that $\frac{d}{dx}\mu(x) = \mu(x)P(x)$ which is simple separable i.e.

$$\frac{1}{\mu(x)}d(\mu(x)) = P(x)dx,$$

$$\int \frac{1}{\mu(x)}d(\mu(x)) = \int P(x)dx,$$

$$\ln[\mu(x)] = \int P(x)dx,$$

$$\mu(x) = e^{\int P(x)dx}.$$

The integrating factor can be obtained provided the function $P(x)$ is integrable.

After $\mu(x)$ has been found, the ODE then reduces to,

$$\frac{d}{dx}\mu(x)y = \mu(x)Q(x).$$

And the solution is;

$$y(x) = \frac{1}{\mu(x)}\int \mu(x)Q(x)\,dx.$$

### (d) Finding the complementary solution and particular integral

In the previous section, ODEs of the form;

$$a_1(x)\frac{dy}{dx} + a_0(x)y = f(x)$$

have been dealt with using the integrating factor method, in general, the integrating factor method works weather $a_1(x)$ and $a_0(x)$ are functions of $x$ or constants.

I now investigate an alternative method to the integrating factor which is in some cases easier to apply, The alternative method also has the advantage that it can be used to solve higher order linear differential equations with constant coefficients; this alternative will only work if $a_1(x)$ and $a_0(x)$ are constants.

Let $a_1(x) = a_1$ and $a_0(x) = a_0$, then

$$a_1 \frac{dy}{dx} + a_0 y = f(x).$$

The first step in this method is to treat the equation as if it were homogeneous, i.e. the RHS of the equation is made zero, $f(x) = 0$.

$$a_1 \frac{dy}{dx} + a_0 y = 0.$$

This is called the **reduced equation** and its solution can be easily evaluated by separation of variables. The solution to the reduced equation is known as the **complementary solution, $y_c$.**

The next step is to try to find, by trial or inspection, a **particular solution**, $y_P$ of the complete equation. **Any** function that satisfies the equation $a_1 \frac{dy}{dx} + a_0 y = f(x)$ is sufficient.

- If $f(x)$ is of the form $ce^{kx}$, where $c$ and $k$ are constants, try a particular integral of the form $y = ae^{kx}$, where $a$ is a constant to be found. This fails when the complementary solution has the same exponential form as the right-hand side of the differential equation. However, it will then be found that the trial function $y = axe^{kx}$ will provide a particular integral.
- If $f(x)$ is of the form $c\cos kx$ or $c\sin kx$, try a PI of the form $y = a\cos kx + b\sin kx$, where $a$ and $b$ are constants to be found.
- If $f(x)$ is a polynomial of degree $n$, try a particular integral is of the form $y = ax^n + bx^{n-1} + \cdots$ where $a, b$, are constants to be found.

Finally, the general solution is the sum of the complementary and particular solutions

$$y = y_c + y_p.$$

The main drawback of this method is that it only works when $f(x)$ takes certain forms i.e. polynomial, exponential and some trigonometric functions.

**(e) Using Special substitutions differential equation**

In many cases, equations can be put into one of the standard forms discussed above (separable, linear, etc.) by a substitution. I will discuss two of the most important general substitutions that are commonly used and their case usage.

i) The Bernoulli equation

This is a differential equation of the form

$$\frac{dy}{dx} + P(x)y = Q(x)y^n \quad \text{where } n \text{ is an integer such that } n > 1$$

Such equation can be converted into an easier form using the substitution

$$z = y^{1-n},$$

$$\frac{dz}{dx} = (1-n)y^{-n}\frac{dy}{dx}.$$

This gives

$$y^{-n}\frac{dy}{dx} + P(x)y^{1-n} = Q(x),$$

$$\frac{1}{1-n}\frac{dz}{dx} + P(x)z = Q(x),$$

$$\frac{dz}{dx} + (1-n)P(x)z = (1-n)Q(x).$$

Now this is linear in z and can be solved by methods described earlier

ii) Homogeneous Equations

*Definition*: (Homogeneous function of degree n)

A function $F(x, y)$ is called homogeneous of degree n if $F(\lambda x, \lambda y) = \lambda^n F(x, y)$. For a polynomial, homogeneous says that all of the terms have the same degree.

Consider $M(x,y)\,dx + N(x,y)\,dy = 0$. Suppose M and N is both homogeneous and of the same degree. Then

$$\frac{dy}{dx} = -\frac{M(x,y)}{N(x,y)},$$

and this suggests that the substitution $y = vx$ might be useful. This substitution once made always changes the unknown function in the differential equation from $y$ to $v$ and makes the differential equation easier to solve by separating variables.

With the substitution $y = vx$, we have

$$-\frac{M(x,y)}{N(x,y)} = R\left(\frac{y}{x}\right),$$

$$\frac{dy}{dx} = R\left(\frac{y}{x}\right) = R(v).$$

since

$$\frac{dy}{dx} = v + x\frac{dv}{dx},$$

we have

$$v + x\frac{dv}{dx} = R(v),$$

$$x\frac{dv}{dx} = R(v) - v.$$

And by separating variables, we obtain

$$\frac{dv}{R(v) - v} = \frac{1}{x}dx.$$

### 3.2.1 Second order linear ODEs with constant coefficients

In this section, the researcher will discuss methods of solving second order linear ordinary differential equations, both homogeneous and non-homogeneous i.e. equations of the type

$$a\frac{d^2y}{dx^2} + b\frac{dy}{dx} + cy = 0,$$

and

$$a\frac{d^2y}{dx^2} + b\frac{dy}{dx} + cy = f(x).$$

In each case, $a$, $b$ and $c$ are given constants and $a \neq 0$ otherwise the differential equation would be of first order.

a) **Second order equation of the form** $a\frac{d^2y}{dx^2} + b\frac{dy}{dx} + cy = 0$

Here I discuss a method of solving a *second order linear homogeneous equation with constant coefficients* i.e. the equation

$$a\frac{d^2y}{dx^2} + b\frac{dy}{dx} + cy = 0,$$

where $a$, $b$ and $c$ are given constants and $a \neq 0$ otherwise the differential equation would be of first order. The form of the general solution of differential equations of this type can be derived by a method which reduces the problem to one of solving first order linear differential equations.

Consider the quadratic equation

$$am^2 + bm + c = 0$$

the coefficients of which are the same as those in the differential equation above. This quadratic is called the **auxiliary equation**. Let the roots be $m_1$ and $m_2$. Then

$$m_1 + m_2 = -\frac{b}{a} \quad \text{and} \quad m_1 m_2 = \frac{c}{a}.$$

Let us set

$$u = \frac{dy}{dx} - m_2 y,$$

then

$$\frac{du}{dx} - m_1 u = \frac{d^2 y}{dx^2} - m_2 \frac{dy}{dx} - m_1 \left(\frac{dy}{dx} - m_2 y\right),$$

$$= \frac{d^2 y}{dx^2} - (m_1 + m_2)\frac{dy}{dx} + m_1 m_2 y,$$

$$= \frac{d^2 y}{dx^2} + \frac{b}{a}\frac{dy}{dx} + \frac{c}{a} y,$$

$$= \frac{1}{a}\left[a\frac{d^2 y}{dx^2} + b\frac{dy}{dx} + cy\right] = 0, \quad a \neq 0$$

$$= 0.$$

Hence, the new variable $u$ satisfies the first order differential equation

$$\frac{du}{dx} = m_1 u$$

which has the general solution $u = c e^{m_1 x}$, where $c$ is an arbitrary constant. Substituting this expression into

$$u = \frac{dy}{dx} - m_2 y$$

Gives

$$\frac{dy}{dx} - m_2 y = c e^{m_1 x},$$

which is a first order linear ODE that can be solved using an integrating factor,

$$IF = e^{\int -m_2 dx} = e^{-m_2 x}.$$

Multiplying $\frac{dy}{dx} - m_2 y = ce^{m_1 x}$ by the integrating factor $e^{-m_2 x}$ yields

$$e^{-m_2 x}\frac{dy}{dx} - m_2 y e^{-m_2 x} = ce^{(m_1 - m_2)x}.$$

The LHS can be reduced to an exact differential

$$\frac{d}{dx}(e^{-m_2 x} y) = ce^{(m_1 - m_2)x}.$$

*Case 1 : For real and distinct roots $m_1 \neq m_2$*

$$e^{-m_2 x} y = \frac{c}{m_1 - m_2} e^{(m_1 - m_2)x} + B$$

Writing $\frac{c}{m_1 - m_2} = A$, *a constant*,

$$y = Ae^{m_1 x} + Be^{m_2 x}.$$

This is therefore, the form of the general solution when the auxiliary equation has real and distinct roots.

*Case 2: For real and equal roots $m_1 = m_2 = m$*

The equation, $\frac{d}{dx}(e^{-m_2 x} y) = ce^{(m_1 - m_2)x}$ becomes

$$\frac{d}{dx}(e^{-mx} y) = c.$$

Hence

$$y = (cx + D)e^{mx},$$

*where c and D are arbitary constants.*

*Case 3: For complex roots $m_1 = \alpha + i\beta, \ m_2 = \alpha - i\beta$*

Considering a case where the roots are complex i.e. $m_1 = \alpha + i\beta$, $m_2 = \alpha - i\beta$, $\alpha$ and $\beta$ are real and $\beta \neq 0$. The general solution is,

$$y = Ae^{(\alpha+i\beta)x} + Be^{(\alpha-i\beta)x},$$

$$y = e^{\alpha x}\left(Ae^{i\beta x} + Be^{-i\beta x}\right).$$

Using Euler's identity,

$$e^{i\theta} = \cos\theta + i\sin\theta,$$

$$y = e^{\alpha x}\left(A(\cos\beta x + i\sin\beta x) + B(\cos-\beta x + i\sin-\beta x)\right),$$

$$y = e^{\alpha x}(C\cos\beta x + D\sin\beta x).$$

where C and D are arbitrary constants.

It is worth noting that when solving equations of the type $a\frac{d^2y}{dx^2} + b\frac{dy}{dx} + cy = 0$ discussed above, we first create an auxiliary equation $am^2 + bm + c = 0$, then find its roots and substitute them into a general solution depending on the nature of the roots.

Note: if one of the linearly independent solutions of $a\frac{d^2y}{dx^2} + b\frac{dy}{dx} + cy = 0$ is known, then the second linearly independent solution can be determined by Abel's theorem or order reduction discussed in section 3.2.2.

### b) Second order equations of the form $\quad a\frac{d^2y}{dx^2} + b\frac{dy}{dx} + cy = f(x)$

In the previous subsection, a method of solving a second order linear homogeneous equation with constant coefficients $a\frac{d^2y}{dx^2} + b\frac{dy}{dx} + cy = 0$ was discussed; I now turn to **a *second order linear non- homogeneous equation with constant coefficients***, $a\frac{d^2y}{dx^2} + b\frac{dy}{dx} + cy = f(x)$ *where* a, b and c are constants.

In section 3.2.0 (d), a method of finding the complementary and particular solutions was used to solve the equation

$$a_1 \frac{dy}{dx} + a_0 y = f(x).$$

The same method can be used to solve second order differential equations of the form

$$a \frac{d^2 y}{dx^2} + b \frac{dy}{dx} + cy = f(x).$$

This method is generally known as the **method of undetermined coefficients**

The procedure is as follows

- Find the solution to the reduced equation $a \frac{d^2 y}{dx^2} + b \frac{dy}{dx} + cy = 0$, this solution is the complementary solution $y_c$. Note that this type of equations have been already dealt with in subsection 3.2.1(a).
- Find a particular solution of the complete equation $a \frac{d^2 y}{dx^2} + b \frac{dy}{dx} + cy = f(x)$, this is the particular solution $y_P$. The particular solution is found by trial and inspection depending on the nature of $f(x)$.
- The general solution is then evaluated as $y = y_c + y_p$.

The method of undetermined coefficients only works if the function $f(x)$ takes up a few special forms. A more general method that works for all forms of the function $f(x)$, known as **variation of parameters** will be discussed in section 3.2.2, this method is advantageous since it works for all forms of $f(x)$, additionally, it works for both second order linear ODEs with constant and variable coefficients

### 3.2.2 Second order linear ODEs with variable coefficients

Consider the differential equation

$$\frac{d^2 y}{dx^2} + P(x) \frac{dy}{dx} + Q(x) y = R(x)$$

where P, Q, and R are, in general, functions of $x$

In a special case where P and Q are constants, then the equation can be solved by the methods described in section 3.2.1(b), but if one or both P and Q are functions of x, then the solution can be found by the Laplace transform method or the power series method.

(a) Variation of parameters

The idea of variation of parameters is to look for functions $v_1(x)$ and $v_2(x)$ such that $v_1 y_1 + v_2 y_2$ satisfy the equation $\frac{d^2 y}{dx^2} + P(x)\frac{dy}{dx} + Q(x)y = R(x)$. Let $y_1$ and $y_2$. If we want $y$ to be a solution to $\frac{d^2 y}{dx^2} + P(x)\frac{dy}{dx} + Q(x)y = R(x)$, then substitution into the D.E will give one condition on $v_1$ and $v_2$.

$$y = v_1 y_1 + v_2 y_2.$$

Differentiating the above equation with respect to x, we obtain,

$$y' = v_1 y_1' + v_2 y_2' + v_1' y_1 + v_2' y_2 = v_1 y_1' + v_2 y_2'$$

since

$$[v_1' y_1 + v_2' y_2 = 0],$$

we find

$$y'' = v_1 y_1'' + v_2 y_2'' + v_1' y_1' + v_2' y_2'$$

Substituting $y'$ and $y''$ into $y'' + P(x)y' + Q(x)y = R(x)$, we obtain

$$\begin{pmatrix} v_1 y_1'' + v_2 y_2'' + v_1' y_1' + v_2' y_2' + P(x)v_1 y_1' \\ + P(x)v_2 y_2' + Q(x)v_1 y_1 + Q(x)v_2 y_2 \end{pmatrix} = R(x)$$

Which simplifies to

$$\begin{pmatrix} v_1(y_1'' + P(x)y_1' + Q(x)y_1) + v_2(y_2'' + P(x)y_2' + Q(x)y_2) \\ + v_1' y_1' + v_2' y_2' \end{pmatrix} = R(x)$$

Since $y_1$ and $y_2$ are solutions, the first two terms are equal to zero, therefore,

$$v_1'y_1' + v_2'y_2' = R(x)$$

$$v_1'y_1 + v_2'y_2 = 0$$

Solving the two equations simultaneously, it follows that

$$v_1' = -\frac{y_2 R(x)}{W}, \qquad v_2' = -\frac{y_1 R(x)}{W},$$

where

$$W = \begin{vmatrix} y_1 & y_2 \\ y_1' & y_2' \end{vmatrix} \neq 0.$$

The Wronskian W can never be zero in this case since $y_1$ and $y_2$ are linearly independent.

The main advantage of this method is that it can be used to solve second order linear differential equations both with constant coefficients or variable coefficients; it can also be used to solve higher order differential equations.

If P and Q are functions of x but $R(x) = 0$, the equation can be solved faster *using Abel's formula* or the *method of order reduction*

(b) Method of order reduction

Let $y_1$ and $y_2$ be linearly independent solutions to the second order linear differential equation and suppose $y_1$ is known.

Consider a subinterval $I$ on which $y_1(x)$ is never zero. So $y_1(x) > 0$ or $y_2(x) < 0$ throughout $I$. We have

$$W = \begin{vmatrix} y_1 & y_2 \\ y_1' & y_2' \end{vmatrix}$$

from Abel's formula;

$W = Ae^{-\int P_1(x)dx}$, for some constant A. (Different A's correspond to different $y_2$'s )

suppose we choose $A = 1$; then

$$y_1 y_2' - y_2 y_1' = e^{-\int P_1(x)dx}$$

$$y_2' - y_2 \frac{y_1'}{y_1} = \frac{e^{-\int P_1(x)dx}}{y_1}$$

Since $y_1$ is known and $y_1 \neq 0$ throughout $I$ the equation above is first order linear in $y_2$, the integrating factor is

$$e^{\int -\frac{y_1'}{y_1}} = e^{-\ln(|y_1|)+C} = \frac{k}{|y_1|}$$

where $k = \pm 1$ depending on the sign of $y_1$, let's take $k = +1$ i.e. use $\frac{1}{y_1}$ as the integrating factor, then

$$\frac{y_2'}{y_1} - y_2 \frac{y_1'}{y_1^2} = \frac{e^{-\int P_1(x)dx}}{y_1^2} .$$

Which simplifies to

$$\frac{y_2}{y_1} = \int \frac{e^{-\int P_1(x)dx}}{y_1^2} dx .$$

Hence

$$y_2 = y_1 \int \frac{e^{-\int P_1(x)dx}}{y_1^2} dx .$$

This method of order reduction can generally be used to reduce an $n^{th}$ order differential equation to an $(n-1)^{th}$ order differential equation.

(c) Laplace Transform method

By applying the Laplace transform, one can change an ordinary differential equation into an algebraic equation, as algebraic equations are generally easier to deal with.

Laplace transforms are used to reduce a differential equation to a simple equation in $s$-space and a system of differential equations to a system of linear equations. We discover that in the case of zero initial conditions, we can solve the system by multiplying the Laplace transform of the input function by the transfer function of the system. Furthermore, if we are only interested in the steady state response, and we have a periodic input, we can find the response by simply adding the response to each of the frequency components of the input signal

From the property of the Laplace transform;

$$\mathcal{L}\{f^n\} = s^n \mathcal{L}\{f\} - \sum_{i=1}^{n} s^{n-1} f^{i-1}(0).$$

Considering the differential equation

$$\sum_{i=0}^{n} a_i f^i(x) = \phi(x),$$

with initial conditions

$$f^i(0) = c_i.$$

Using the linearity of the Laplace transform it is equivalent to rewrite the equation $\sum_{i=0}^{n} a_i f^i(x) = \phi(x)$, as

$$\sum_{i=0}^{n} a_i \mathcal{L}\{f^i(x)\} = \mathcal{L}\{\phi(x)\},$$

which gives

$$\mathcal{L}\{f(x)\} \sum_{i=0}^{n} a_i s^i - \sum_{i=1}^{n} \sum_{j=1}^{i} a_i s^{i-j}(0) = \mathcal{L}\{\phi(x)\}.$$

Solving the equation for $\mathcal{L}\{f(t)\}$ and substituting $f^i(0) = c_i$

$$\mathcal{L}\{f(x)\} = \frac{\mathcal{L}\{\phi(x)\} \sum_{i=1}^{n} \sum_{j=1}^{i} a_i s^{i-j} c_{j-1}}{\sum_{i=0}^{n} a_i s^i}$$

The solution for $f(t)$ is obtained by applying the inverse Laplace transform to $\mathcal{L}\{f(t)\}$

Since all the initial conditions are zero,

$$f^i(0) = c_i = 0 \qquad \forall_i \in \{0,1,2, \ldots, n\}.$$

Then

$$f(x) = \mathcal{L}^{-1} \left\{ \frac{\mathcal{L}\{\phi(x)\}}{\sum_{i=0}^{n} a_i s^i} \right\}.$$

(d) Power Series method

Linear second order homogeneous ODEs with polynomials as functions can often be solved by expanding functions around ordinary or specific points. In general, such a solution assumes a power series with unknown coefficients, then substitutes that solution into the differential equation to find a recurrence relation for the coefficients.

Consider the second order linear differential equation,

$$a_2(x) \frac{d^2 y}{dx^2} + a_1(x) \frac{dy}{dx} + a_0(x) y = 0,$$

since

$$a_2(x) \neq 0,$$

it follows that

$$\frac{d^2y}{dx^2} + \frac{a_1(x)}{a_2(x)}\frac{dy}{dx} + \frac{a_0(x)}{a_2(x)}y = 0 \ .$$

Suppose $\frac{a_1(x)}{a_2(x)}$ and $\frac{a_0(x)}{a_2(x)}$ are analytic functions, using the power series method, we can form a power series solution,

$$y(x) = \sum_{n=0}^{\infty} a_n x^n.$$

If for some value of x, $a_2(x) = 0$, then a variation of the power series method called the Frobenius method is used.

The steps for using the power series method are as follows **[35]**;

i)      Assume the differential equation has a solution of the form

$$y(x) = \sum_{n=0}^{\infty} a_n x^n.$$

ii)      Differentiate the power series term by term to get

$$y'(x) = \sum_{n=0}^{\infty} n a_n x^{n-1}$$

and

$$y''(x) = \sum_{n=0}^{\infty} n(n-1) a_n x^{n-2} \ .$$

iii)      Substitute the power series expressions into the differential equation.

iv)      Re-index sums as necessary to combine terms and simplify the expression.

v)      Equate coefficients of like powers of *x* to determine values for the coefficients *an* in the power series.

vi)      Substitute the coefficients back into the power series and write the solution.

## 3.3 Numerical methods of solving an ordinary differential equation

Numerical methods are those used to find numerical approximations to the solutions of ordinary differential equations (ODEs). Their use is also known as "numerical integration", although this term is sometimes taken to mean the computation of integrals.

Although many of the differential equations which result from modeling real-world problems can be solved analytically, there are many others which cannot. In general, when the modeling leads to a **linear** differential equation, the prospects of obtaining an exact mathematical solution are good. However, non-linear differential equations present much greater difficulty and exact solutions can seldom be obtained. There is a need, therefore, for numerical methods that can provide approximate solutions to problems which would otherwise be intractable. The advent of powerful computers capable of performing calculations at very high speed has led to a rapid development in this area and there are now many numerical methods available.

Thus numerical methods are useful because many differential equations in practice cannot be solved using symbolic computation ("analysis") and so for practical purposes, such as in engineering – a numeric approximation to the solution is often sufficient. The algorithms studied in this research can be used to compute such an approximation.

Consider the ODE below

$$\frac{dy}{dx} = f(x, y)$$

Subject to initial conditions $y(x_0) = y_0$ and suppose the solution $y(x)$ for $x \geq x_0$ is required.

Suppose the graph shown below represents the exact solution

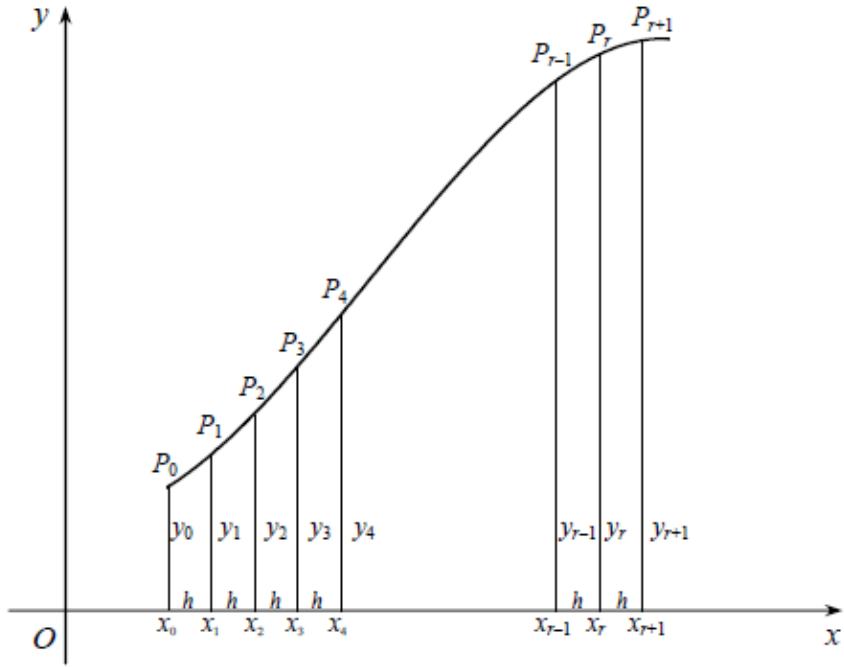

Figure 1

Let $x_0, x_1, x_2, \ldots$ be points along the x-axis which are equally spaced a distance $h$ apart and their corresponding y-values $y(x_0), y(x_1), y(x_2), \ldots$ be represented as $y_0, y_1, y_2, \ldots$ and let $P_0, P_1, P_2, \ldots$ be points on the curve with coordinates $(x_0, y_0), (x_1, y_1), (x_2, y_2) \ldots$ respectively. Therefore, in general, $y(x_r) = y_r$ where $x_r = x_0 + rh$ and $P_r$ is a point with coordinate $(x_r, y_r)$.

The basic idea of all the numerical methods is to use a step-by-step procedure to obtain **approximations** to the values of $y_0, y_1, y_2, \ldots$ successively. The interval, $h$, between successive x-values is called the **step length**. Some of the methods of obtaining these approximations are suggested by geometrical considerations of the graph.

### 3.3.1 Euler's method

Consider figure 1 above, one simple way of obtaining an approximation to the value of $y_1$ is to assume that the part of the curve between $P_0$ and $P_1$ is a straight-line segment with gradient equal to the gradient of the curve at $P_0$.

Since $\frac{dy}{dx} = f(x, y)$, the gradient of the curve at $P_0$ is $f(x_0, y_0)$. Hence with this approximation,

$$\frac{y_1 - y_0}{h} = f(x_0, y_0),$$

giving

$$y_1 = y_0 + hf(x_0, y_0).$$

Using this approximation, we obtain the value of $y$ at $P_1$, The process can be repeated, assuming that the part of the curve between $P_1$ and $P_2$ is a straight line segment with gradient equal to the value of $\frac{dy}{dx}$ at $P_1$. This gives

$$y_2 = y_1 + hf(x_1, y_1).$$

More generally,

$$y_{r+1} = y_r + hf(x_r, y_r), \quad r = 0,1,2,\ldots.$$

This is *Euler's formula*. Successive calculation of values of $y$ using this formula is known as *Euler's method*. It is clear from the nature of the linear approximation on which this method is based that the step length $h$ needs to be fairly small to achieve reasonable accuracy.

### 3.3.2 The improved Euler formula

The Euler formula derived above is based on the assumption that the gradient of the line segment joining any two successive points, $P_r$ and $P_{r+1}$, on the graph of $y(x)$ is equal to the value of $\frac{dy}{dx}$ at $P_r$. Another formula, which is considerably more accurate than Euler's, can be obtained by assuming that the gradient of $P_r P_{r+1}$ is equal to the average of the values of $\frac{dy}{dx}$ at $P_r$ and $P_{r+1}$. With this assumption,

$$\frac{y_{r+1} - y_r}{h} = \frac{1}{2}\{f(x_r, y_r) + f(x_{r+1}, y_{r+1})\}.$$

Hence,

$$y_{r+1} = y_r + \frac{h}{2}\{f(x_r, y_r) + f(x_{r+1}, y_{r+1})\}.$$

This formula cannot be used directly to calculate $y_{r+1}$ because $y_{r+1}$ is needed in order to evaluate $f(x_{r+1}, y_{r+1})$ on the right-hand side. To overcome this problem, $y_{r+1}$ on the right-hand side is replaced by a first estimate, denoted by $y_{r+1}^*$, giving

$$y_{r+1} = y_r + \frac{h}{2}\{f(x_r, y_r) + f(x_{r+1}, y_{r+1}^*)\}.$$

There are a number of ways of obtaining a value for $y_{r+1}^*$. The simplest is to use Euler's formula, which gives

$$y_{r+1}^* = y_r + hf(x_r, y_r),$$

### 3.3.3 Error Analysis

In the previous subsections, I have derived the Euler and Improved Euler method by geometrical constructions, I will now derive them analytically, this alternative way of derivation offers an insight about the errors involved.

The maclaurin series for $y(h)$ is

$$y(h) = y(0) + hy'(0) + \frac{h^2}{2!}y''(0) + \frac{h^3}{3!}y'''(0) + \cdots$$

If the origin is transferred to point $x_r$, this becomes

$$y(h + x_r) = y(x_r) + hy'(x_r) + \frac{h^2}{2!}y''(x_r) + \frac{h^3}{3!}y'''(x_r) + \cdots$$

which is a taylor series expansion

if $\qquad y'(x_r) = f(x_r, y_r)$

the taylor series may be expressed as

$$y_{r+1} = y_r + hf(x_r, y_r) + \frac{h^2}{2!}y''(x_r) + \frac{h^3}{3!}y'''(x_r) + \cdots$$

For a sufficiently small choice of $h$, the terms $h^2$, $h^3$, ...can be neglected and the equation reduces to $y_{r+1} = y_r + hf(x_r, y_r)$, which is the Euler's formula.

The error incurred by neglecting terms in an expansion is called the **truncation error,** The above calculation suggests that a smaller value of step length, h will yield a low truncation error, this will be shown in chapter 4 by implementing some case studies with MATLAB and investigating the results for different value of h.

### 3.3.4 Runge-Kutta Methods

From Euler's formula, $y_{r+1} = y_r + hf(x_r, y_r)$,

We denote $k_1 = f(x_r, y_r)$

Therefore, going forth to the approximation to $x_r + \frac{1}{2}h$

$$k_2 = f(x_r + \frac{1}{2}h, \ y_r + \frac{1}{2}hk_1),$$

$$k_3 = f(x_r + \frac{1}{2}h, \ y_r + \frac{1}{2}hk_2),$$

$$k_4 = f(x_r + \frac{1}{2}h, \ y_r + \frac{1}{2}hk_3),$$

hence

$$y_{r+1} = y_r + \frac{h}{6}\{k_1 + 2k_2 + 2k_3 + 3k_4\}.$$

This gives a local truncation error proportional to $h^5$ and since $h < 1$ in most cases, it means the error is much less with the fourth order Runge-Kutta method.

Each of the $k_i$, $i = 1, 2, ...$ represents a k[th] order Runge-Kutta method, the fourth order is the most stable and easy to implement

### 3.3.5 MATLAB Implementation of numerical Algorithms

The MATLAB scientific computing software has a number of inbuilt 'functions' for solving ordinary differential equations to a reasonable degree of accuracy, these include ode45, ode23, ode115 and many others.

The researcher however did not use these inbuilt 'ode solvers' but rather implemented the Euler and Runge Kutta formulas from scratch using MATLAB programming Language. This implementation allowed the researcher to study the effect of different step sizes on the accuracy of the solution obtained, it also enabled visual representation of results through graphs. The source code for the programs that were written to generate numerical solutions are attached at the appendix of this document.

# CHAPTER FOUR

# REPRESENTATION AND ANALYSIS OF RESULTS

## 4.0 Introduction

In this chapter, we implement numerical and analytical methods of solving ordinary differential equations to specific problems. Here, the researcher undertook a case study of a numerical problem, computed its solution using analytical methods (whenever possible) to find the exact solution and using numerical methods to find the approximate solution, the error bounds between different numerical methods and exact values for selected step sizes were computed to analyze their effect on the accuracy.

MATLAB mathematical computing software was used to compute numerical solutions and plot their graphs to visualize their meaning and deviation between the solutions.

## 4.1 Case study

**Problem**: Consider the differential equation $y'(x) = -6y, \ given \ y(0) = 1$

The exact solution to this problem is obtained analytically by separation of variables as $\boldsymbol{y = e^{-6x}}$

$$\frac{dy}{dx} = -6y$$

By separating variables, we get;

$$\frac{dy}{y} = -6dx,$$

On integrating both sides, we obtain;

$$\int \frac{dy}{y} = -6 \int dx$$

$$\ln y = -6x + k$$

$$y = e^{-6x+k}.$$

Applying the given conditions,

$y(0) = 1$

$$1 = e^k,$$

thus

$$k = 0.$$

And

$$y = e^{-6x}$$

is the exact solution to the differential equation at the given initial conditions.

In this case study, I have used the method of separating variables since the algebra allows for 'like terms' to be gathered to one side. It is also possible to use other analytical methods although in this case, the above method is fastest.

The equation was then solved using Euler and 4th order Runge-Kutta methods for different step sizes and the errors in each case the error in each case calculated. The table below shows the results obtained.

*Table 4.1: Comparison between exact value and that obtained by numerical methods, Euler method and runge kutta method for step size, h=0.1*

| $x_n$ | Exact Value (analytical solution) $y = e^{-6x_n}$ | Euler's Method, h=0.1 $y(x_n)$ | Absolute error | 4$^{th}$ Runge Kutta Method, h=0.1 $y(x_n)$ | Absolute error |
|---|---|---|---|---|---|
| 0.0 | 1.000000000000000 | 1.000000000000000 | 0.000000000 | 1.000000000000000 | 0.0000000000 |
| 0.1 | 0.548811636094026 | 0.400000000000000 | 0.148811636 | 0.549400000000000 | 0.000588364 |
| 0.2 | 0.301194211912202 | 0.160000000000000 | 0.141194212 | 0.301840360000000 | 0.000646148 |
| 0.3 | 0.165298888221586 | 0.064000000000000 | 0.101298888 | 0.165831093784000 | 0.000532206 |
| 0.4 | 0.090717953289412 | 0.025600000000000 | 0.065117953 | 0.091107602924930 | 0.00038965 |
| 0.5 | 0.049787068367864 | 0.010240000000000 | 0.039547068 | 0.050054517046956 | 0.000267449 |
| 0.6 | 0.027323722447293 | 0.004096000000000 | 0.023227722 | 0.027499951665598 | 0.000176229 |
| 0.7 | 0.014995576820477 | 0.001638400000000 | 0.013357177 | 0.015108473445079 | 0.000112897 |
| 0.8 | 0.008229747049020 | 0.000655360000000 | 0.007574387 | 0.008300595310727 | 0.000070848 |
| 0.9 | 0.004516580942613 | 0.000262144000000 | 0.004254437 | 0.004560347063713 | 0.000043766 |
| 1.0 | 0.002478752176666 | 0.000104857600000 | 0.002373895 | 0.002505454676804 | 0.000026702 |

*Table 4.2: Comparison between exact value and that obtained by numerical methods, Euler method and runge kutta method for step size, h=0.05*

| $x_n$ | Exact Value (analytical solution) $y = e^{-6x_n}$ | Euler's Method, h=0.05 $y(x_n)$ | Absolute error | 4th Runge Kutta Method, h=0.05 $y(x_n)$ | Absolute error |
|---|---|---|---|---|---|
| 0.0 | 1.000000000000000 | 1.000000000000000 | 0.000000000 | 1.000000000000000 | 0.00000000 |
| 0.1 | 0.548811636094026 | 0.490000000000000 | 0.058811636 | 0.548840201406250 | $2.85653 \times 10^{-5}$ |
| 0.2 | 0.301194211912202 | 0.240100000000000 | 0.061094212 | 0.301225566679653 | $3.13548 \times 10^{-5}$ |
| 0.3 | 0.165298888221586 | 0.117649000000000 | 0.047649888 | 0.165324700685173 | $2.58125 \times 10^{-5}$ |
| 0.4 | 0.090717953289412 | 0.057648010000000 | 0.033069943 | 0.090736842021478 | $1.88887 \times 10^{-5}$ |
| 0.5 | 0.049787068367864 | 0.028247524900000 | 0.021539543 | 0.049800026650035 | $1.29583 \times 10^{-5}$ |
| 0.6 | 0.027323722447293 | 0.013841287201000 | 0.013482435 | 0.027332256656642 | $8.53421 \times 10^{-6}$ |
| 0.7 | 0.014995576820477 | 0.006782230728490 | 0.008213346 | 0.015001041248319 | $5.46443 \times 10^{-6}$ |
| 0.8 | 0.008229747049020 | 0.003323293056960 | 0.004906454 | 0.008233174500031 | $3.42745 \times 10^{-6}$ |
| 0.9 | 0.004516580942613 | 0.001628413597910 | 0.002888167 | 0.004518697150810 | $2.11621 \times 10^{-6}$ |
| 1.0 | 0.002478752176666 | 0.000797922662976 | 0.00168083 | 0.002480042654344 | $1.29048 \times 10^{-6}$ |

*Table 4.3: Comparison between exact value and that obtained by numerical methods, Euler method and runge kutta method for step size, h=0.01*

| $x_n$ | Exact Value (analytical solution) $y = e^{-6x_n}$ | Euler's Method, h=0.01 $y(x_n)$ | error | 4th Runge Kutta Method, h=0.01 $y(x_n)$ | error |
|---|---|---|---|---|---|
| 0.0 | 1.000000000000000 | 1.000000000000000 | 0.000000000 | 1.000000000000000 | 0.000000000 |
| 0.1 | 0.548811636094026 | 0.538615114094900 | 0.010196522 | 0.548811673481706 | $3.73877 \times 10^{-8}$ |
| 0.2 | 0.301194211912202 | 0.290106241131462 | 0.011087971 | 0.301194252949791 | $4.10376 \times 10^{-8}$ |
| 0.3 | 0.165298888221586 | 0.156255606166665 | 0.009043282 | 0.165298922004447 | $3.37829 \times 10^{-8}$ |
| 0.4 | 0.090717953289412 | 0.084161631143426 | 0.006556322 | 0.090717978009983 | $2.47206 \times 10^{-8}$ |
| 0.5 | 0.049787068367864 | 0.045330726560729 | 0.004456342 | 0.049787085326535 | $1.69587 \times 10^{-8}$ |
| 0.6 | 0.027323722447293 | 0.024415814458512 | 0.002907908 | 0.027323733615832 | $1.11685 \times 10^{-8}$ |
| 0.7 | 0.014995576820477 | 0.013150726690291 | 0.00184485 | 0.014995583971473 | $7.15100 \times 10^{-9}$ |
| 0.8 | 0.008229747049020 | 0.007083180156722 | 0.001146567 | 0.008229751534220 | $4.4852 \times 10^{-9}$ |
| 0.9 | 0.004516580942613 | 0.003815107888268 | 0.000701473 | 0.004516583711834 | $2.76922 \times 10^{-9}$ |
| 1.0 | 0.002478752176666 | 0.002054874770524 | 0.000423877 | 0.002478753865312 | $1.68865 \times 10^{-9}$ |

**Fig. 2. A graph showing the comparison between the Exact solution and the solution obtained using Euler's method for different step sizes**

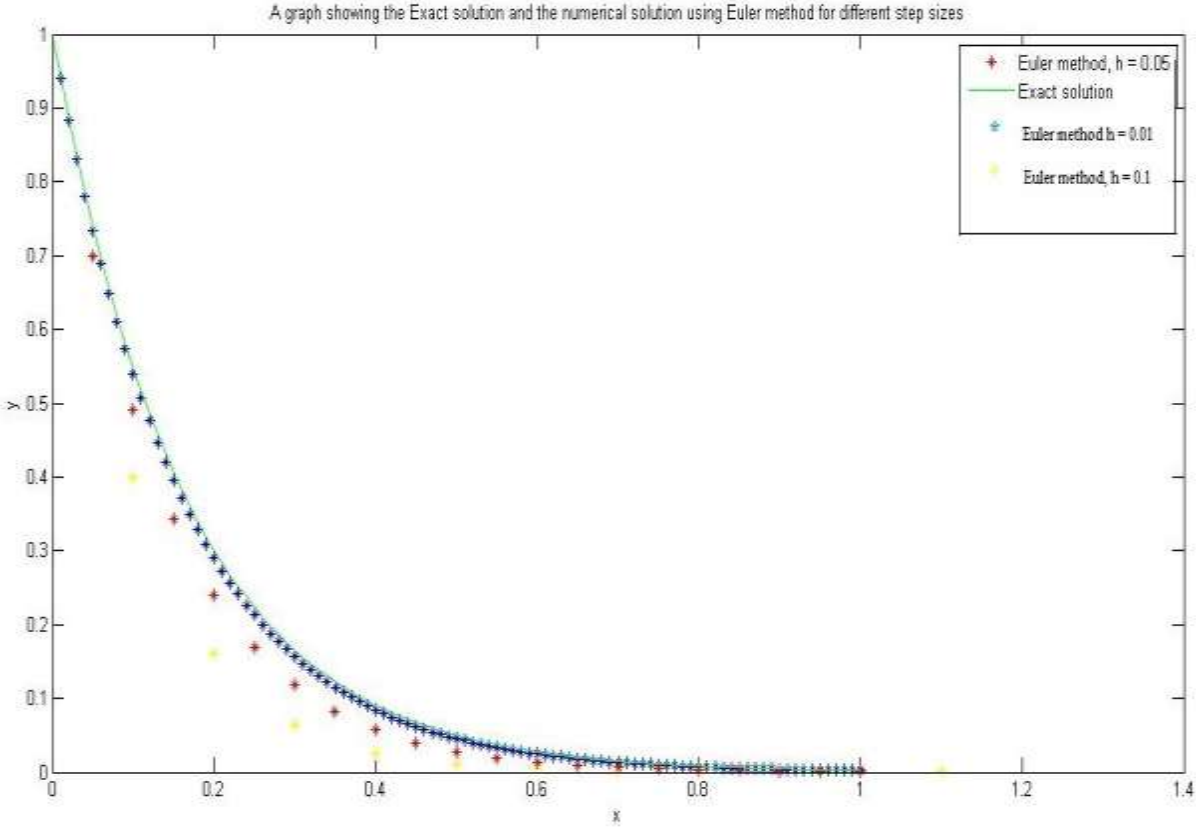

**Fig. 3. A graph showing the comparison between the Exact solution and the solution obtained using Runge Kutta method for different step sizes**

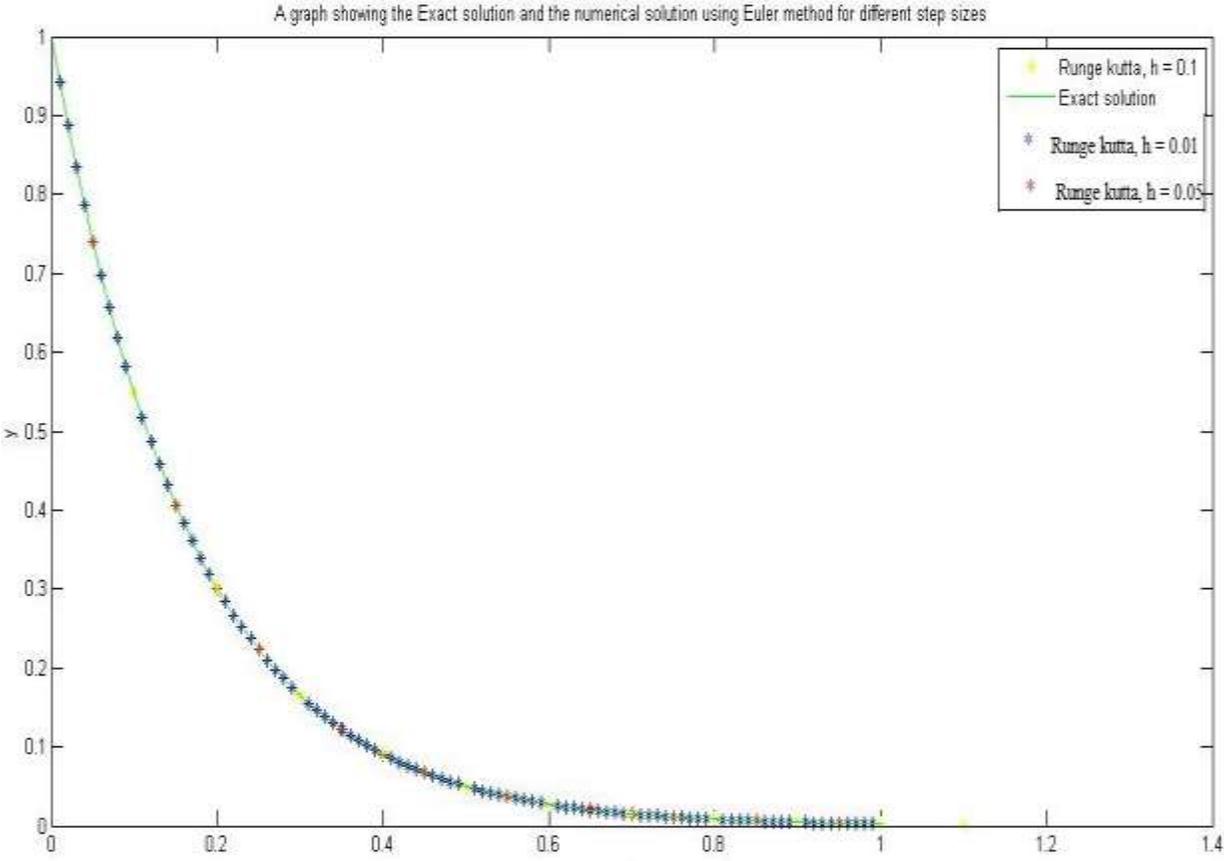

### 4.1.1 Discussion of Results

From the tables 4.1, 4.2 and 4.3 above, the approximate solution is calculated using Euler and Runge Kutta methods with step sizes 0.1, 0.05 and 0.01 respectively and maximum errors also are calculated at specified step sizes. From the tables, for each method, it is seen that a numerical solution converges to the exact solution if the step size leads to decreased errors such that in the limit when the step size approaches to zero the errors go to zero. Also, we see that the Runge Kutta approximations for same step size converge firstly to exact solution than the Euler approximations.

The smaller the step size, the more the numerical solution approaches the exact solution on the graph. This shows that the small step size provides the better approximation.

Finally, we observe that the fourth order Runge Kutta method is converging faster than the Euler method and it is the most effective method for solving initial value problems for ordinary differential equations.

# CHAPTER FIVE

# CONCLUSION AND RECOMMENDATIONS

## 5.0 Introduction

This section involves making conclusions and recommendations basing on the findings of the research study carried out.

## 5.1 Conclusion

The main purpose of this research study was to bring together the different analytical and numerical methods for solving ordinary differential equations. Emphasis was laid on first and second order equations although some of the methods discussed primarily apply even to higher order differential equations.

Scientific computation was used in this study to ease work in doing repetitive calculations and avoid human errors. Euler and Runge kutta algorithms were implemented using MATLAB

## 5.2 Accuracy of the Euler and Runge Kutta methods.

The two main numerical methods studied in this research were the Euler method and the $4^{th}$ Runge kutta method. Both methods are useful in obtaining approximate solutions to differential equations, in fact, they become even more useful in cases where the analytical solution does not exist. It is therefore appropriate to conclude with a note on the accuracy of these methods.

When both methods were used to solve the same initial value problems, it was found that for the same step size, h, the error associated with Euler methods was quite larger than that obtained from the $4^{th}$ Runge Kutta method which agrees with the theory.

This suggests that the Runge kutta method is far better than the single step Euler method. Infact for step sizes less than h=0.05, using the Runge kutta method, the graph for numerical solution was almost superimposed on the exact analytic solution.

For both methods, it was observed that better results are obtained when the step size is made smaller, although this increases the number of calculation cycles to be made and makes the computation a bit intensive.

## 5.3. Recommendations for future researchers

From the survey results of this research and the problems encountered during the study, the researcher recommends the following to future researchers;

1) Future researchers should use the MAPLE scientific computation software to find analytic solutions to ordinary differential equations. Although this software may not give analytic solutions to all differential equations, it saves time in finding solutions for equations that can be solved analytically.

2) For ordinary differential equations whose solutions cannot be solved analytically, the numerical solutions should be computed using MATLAB. Specifically, by either coding these numerical methods from scratch or using inbuilt 'ode solvers'

3) Future researchers should compare the accuracy of different 'ode solvers' in MATLAB such as ode 45, ode23 and others so that they are able to determine which solver is best for a given kind of equations.

4) Future researchers are urged to study more about 'stiff' differential equations as they often times occur in problems obtained from natura phenomena, science, engineering, economics and other models.

5) Because of the high cost associated with the purchase of scientific computing packages such as MAPLE, WOLFRAM & MATLAB which at times hinders their use, future researchers are encouraged to study more about the implementation of numerical methods using general purpose programming languages such as python, C, C++, R and others. They should also look into the use of open source mathematical packages such as OCTAVE, etc.

6) Future researchers should also give more attention to finding solutions to systems of Ordinary differential equations and higher order differential equations using both the analytical methods, whenever possible and also Numerical methods

# APPENDIX

## MATLAB program for computing the solution using Runge Kutta method (4th order)

```
clear all

clc

f=@(x,y)-6*y; %equation to be solved.

x0=0; %inititial value of x

y0=1;   %inititial value of y

xn=1; % end value of x

h=0.1; % Step length

 fprintf('\n x      y ');

while x0<=xn

   fprintf('\n%4.3f  %18.15f ',x0,y0); %values of x and y

   k1=h*f(x0,y0);

   k2=h*f(x0+h/2,y0+k1/2);

   k3=h*f(x0+h/2,y0+k2/2);

   x1=x0+h;

   k4=h*f(x1,y0+k3);

   y1=y0+(k1+2*(k2+k3)+k4)/6;

   x0=x1;

   y0=y1;
```

```matlab
% Plotting the numerical solution

p1=plot(x0,y0,'y*');

hold on

end

%Exact value calculation and plotting

x=0:0.01:1;

y=exp(-6*x);

p2=plot(x,y,'g-');

title('A graph showing the Exact solution and the numerical solution using Euler method for different step sizes')

xlabel('x')

ylabel('y')

legend([p1 p2],{'Runge kutta, h = 0.1','Exact solution'},'Location','Northeast')

hold on
```

## MATLAB program for computing the solution using Euler's method

```
clear all

clc

f=@(x,y) -6*y; % dy/dx = -6*y

x0=0; % initial value of x

y0=1; % initial value of y

xn=1.00;

h=0.05; % Step size

 %calculation loop

 fprintf('\n x      y ');

while x0<=xn

   fprintf('\n%4.3f  %18.15f ',x0,y0); %values of x and y

   y1=y0+h*f(x0,y0);

   x1=x0+h;

   x0=x1;

   y0=y1;

   p1 = plot(x0,y0,'r*'); % Plot for numerical solution

   hold on

end

x=0:0.01:1.00;

y=exp(-6*x);
```

```
p2= plot(x,y,'g-'); % Plot for exact solution

title('A graph showing the Exact solution and the numerical solution using Euler method for different step sizes')

xlabel('x')

ylabel('y')

legend([p1 p2],{'Euler method, h = 0.05','Exact solution'},'Location','Northeast')

hold on
```